\documentclass[10pt,a4paper,reqno]{amsart}
\usepackage{amsfonts,amsthm,latexsym,amsmath,amssymb,amscd,amsmath}
\usepackage{epsfig,psfrag}
\usepackage{color}

\usepackage{geometry}                		% See geometry.pdf to learn the layout options. There are lots.
\usepackage{graphicx}				% Use pdf, png, jpg, or epsÃÂÃÂ§ with pdflatex; use eps in DVI mode
\usepackage{amssymb}
\usepackage[colorlinks = true]{hyperref}

\newtheorem{Theorem}{Theorem}[section]
\newtheorem{Lemma}[Theorem]{Lemma}

\newtheorem{problem}{Problem}

\newtheorem{theo+}[Theorem]{Theorem}
\newtheorem{prop+}[Theorem]{Proposition}
\newtheorem{coro+}[Theorem]{Corollary}
\newtheorem{lemm+}[Theorem]{Lemma}
\newtheorem{conjecture}[Theorem]{Conjecture}
\theoremstyle{definition}
\newtheorem{defi+}[Theorem]{Definition}

\newtheorem{Rmk}[Theorem]{Remark}

\newcommand \bC {\mathbb C}

\newcommand \al {\alpha}
\newcommand \be {\beta}
\newcommand \la {\lambda}

\newcommand{\bP}{\mathbb P}
\newcommand \rk {\text {rk}}
\newcommand \Z {\mathcal Z}
\newcommand \bN {\mathbb N}

\newcommand \HF {\mathrm{HF}}
\newcommand \HS {\mathrm{HS}}

\begin{document}

          \title[Algebraic Stories from One and from  the Other Pockets]
          {Algebraic Stories from One  and from  the Other Pockets}

\author[R.~Fr\"oberg, S.~Lundqvist]{Ralf Fr\"oberg, Samuel Lundqvist}
\address[R.~Fr\"oberg, S.~Lundqvist, B.~Shapiro]{   Department of Mathematics,
            Stockholm University,
            S-10691, Stockholm, Sweden}
\email{ralff@math.su.se, samuel@math.su.se, shapiro@math.su.se}

\author[A.~Oneto]{Alessandro Oneto}
\address[A.~Oneto]{Department of Mathematics, Universitat Polit\`ecnica de Catalunya, Barcelona, Spain} \email{alessandro.oneto@inria.fr, aless.oneto@gmail.com}

\author[B.~Shapiro]{Boris Shapiro}

\begin{abstract}   In what follows, we present a number of  questions which were posed on the problem solving seminar in  algebra at Stockholm University during the period Fall 2014 -- Spring 2017. %These problems  mainly deal with the Waring problem for forms and the HIlbert series of different ideals. %
\end{abstract}

\dedicatory{To our dear colleague, late Jan-Erik Roos} 

\maketitle

{\let\thefootnote\relax\footnotetext{The title alludes to the famous collection of mystery novels ``Pov\'idky z jedn\'e a z druh\'e kapsy" by K.~\v Capek.}}

\section{The Waring problem for complex-valued forms}
The following famous  result on binary forms was   proven by J.~J.~Sylvester in 1851. (Below we use the terms ``forms" and ``homogeneous polynomials" as synonyms.) 

\begin{Theorem}[Sylvester's Theorem \cite{Syl2}]\label{th:Sylv}
\rm{(i)} A general binary form $f$ of odd degree $k=2s-1$ with complex coefficients can be written as
$$ f(x, y) =\sum_{j=1}^s(\al_jx+\be_jy)^k.$$

\noindent
\rm{(ii)} A general binary form $f$ of even degree $k=2s$ with complex coefficients can be written as 
$$ f(x,y)=\la x^k +\sum_{j=1}^s(\al_jx+\be_jy)^k.$$
\end{Theorem}

%This result motivates the study of $k$-ranks. %following definition. 
Sylvester's result was the starting point of the study of the so-called \emph{Waring problem for polynomials} which we discuss below.

\medskip
 Let $S = \bC[x_1,\ldots,x_n]$ be the polynomial ring  in $n$ variables with complex coefficients. With respect to the standard grading, we have $S = \bigoplus_{d\geq 0} S_d$, where $S_d$ denotes the vector space of  all forms of degree $d$.

\begin{defi+}\label{def1}
Let $f$ be a form of degree $k$ in $S$. A presentation of $f$ as a sum of $k$-th powers of linear forms, i.e., $f=l_1^k+\ldots+l_s^k$, where $l_1,\ldots,l_s \in S_1$, is called a \emph{Waring decomposition} of $f$. 
The minimal length of such a decomposition is called the \emph {Waring rank} of $f$, and we denote it as $\rk(f)$.
 By $\rk^\circ(k,n)$ we denote  the Waring rank of a {\em general} complex-valued form of degree $k$ in $n$ variables.
\end{defi+} 

%$k=2s-1$ equals $s$, while the Waring rank of a general binary form of degree $k=2s$ equals $s+1.$

\begin{Rmk}
Besides being a natural question from the point of view of algebraic geometry, the Waring problem  for polynomials is partly motivated by its celebrated prototype, i.e., the Waring problem for natural numbers.  The latter was posed in 1770 by a British number theorist E. Waring %(1736 -- 1798), t 
who claimed that, for any positive integer $k$, there exists a minimal number $g(k)$ such that every natural number can be written as sum of at most $g(k)$ $k$-th powers of positive integers. The famous Lagrange's four-squares Theorem (1770) claims that $g(2) = 4$  while the existence of $g(k)$, for any integer $k\ge 2$, is due to D. Hilbert (1900). Exact values of $g(k)$ are currently known only in a few cases.
\end{Rmk}

%Let $S = \bC[x_1,\ldots,x_n]$ be the polynomial ring with complex coefficients in $n$ variables. With respect to the standard grading, we have $S = \bigoplus_{d\geq 0} S_d$, where $S_d$ denotes the vector space of homogeneous polynomials, or {\it forms}, of degree $d$.

\subsection{Generic $k$-rank}
In terms of Definition~\ref{def1}, Sylvester's Theorem claims that the Waring rank of a general binary complex-valued  form of degree $k$ equals $\left\lfloor \frac{k}{2} \right\rfloor$. 
More generally, the important result of J. Alexander and A. Hirschowitz \cite{al-hi} completely describes the Waring rank $\rk^\circ(k,n)$ of  general forms of any degree and in any number of variables. 

\begin{Theorem}[Alexander-Hirschowitz Theorem, 1995]\label{th:AH} For all pairs of positive integers $(k,n)$,  the generic Waring rank  $\rk^\circ(k,n)$ is given by 
\begin{equation}\label{rgen}
\rk^\circ(k,n)=\left\lceil \frac{\binom{n+k-1}{n-1}}{n} \right\rceil,
\end{equation}
except for the following cases:
\begin{enumerate}
	\item $k = 2$, where $\rk^{\circ}(2,n) = n$;
	\item $k=4,\; n=3,4,5$; and $k=3, n=5$, where, $r_{gen}(k,n)$ equals the r.h.s of \eqref{rgen} plus $1$.
\end{enumerate}
\end{Theorem} 

Going further,  R.F. and B.S. jointly with G.~Ottaviani considered the following natural version of the Waring problem for complex-valued forms, see \cite{FOS}. 

\begin{defi+}
Let $k,d$ be positive integers. Given a form $f$ of degree $kd$, a {\em $k$-Waring decomposition} is a presentation of $f$ as a sum of $k$-th powers of forms of degree $d$, i.e., $f = g_1^k + \cdots + g_s^k$, with $g_i \in S_d$. The minimal length of such an expression is called the {\em $k$-rank} of $f$ and is denoted by $\rk_k(f)$.
We denote by $\rk_k^\circ(kd,n)$ the $k$-rank of a general complex-valued form of degree $kd$ in $n$ variables.
 \end{defi+}
 
 In this notation, the case $d = 1$ corresponds to the classical Waring rank, i.e., if $k = \deg(f)$, then $\rk(f) = \rk_{k}(f)$ and $\rk^\circ(k,n)=\rk_k^\circ(k,n)$. Since the case $k = 1$ is trivial,  we assume below that $k \geq 2$.

%The number $r_{max}(k,d,n)$ is called the {\it maximal $k$-rank}. 

%can be represented as the sum of at most  $r_{max}(k,d,n)$ $k$-th powers of  forms of degree $d$.

\begin{problem}
Given a triple of positive integers $(k,d,n)$, calculate  $\rk_k^\circ(kd,n).$
\end{problem}

The main  result of \cite{FOS}  states  that, for any triple $(k,d,n)$ as above,   
\begin{equation}\label{ineq:main}
\rk^\circ_k(kd,n)\le k^{n-1}.
\end{equation}
At the same time, by a simple parameter count, one has an obvious  lower bound for  $\rk^\circ(k,n)$ given by
\begin{equation}\label{ineq:lower bound}
\rk^\circ_k(kd,n)\geq \left\lceil \frac{\binom{n+kd-1}{n-1}} {\binom{n+d-1}{n-1}} \right\rceil.
\end{equation}

\medskip
A remarkable fact about the upper bound given by \eqref{ineq:main} is that it is independent of $d$. Therefore, since the right-hand side of \eqref{ineq:lower bound} equals  $k^n$ when $d \gg 0$,  we get that for large values of $d$, the bound in \eqref{ineq:main} is actually sharp. As a consequence of this remark, for any fixed $n\ge 1$ and $k\ge 2$, there exists a positive integer $d_{k,n}$ such that $\rk^\circ_k(kd,n)= k^n,$ for all $d\ge d_{k,n}$. 

In the case of binary forms, it has been proven that \eqref{ineq:lower bound} is actually an equality \cite{Re1, LuOnReSh}. Exact values of $d_{k,n}$, and the behaviour of $\rk_k(kd,n)$ for $d \leq d_{k,n}$, have also been computed  in a few other cases, see \cite[Section 3.3]{On}. These results agree with the following illuminating conjecture suggested by G.~Ottaviani in 2014.

\begin{conjecture}\label{conj:main}
The $k$-rank of a general form of degree $kd$ in $n$ variables is given by 
\begin{equation}\label{eq:main}
\rk_k^\circ(kd,n)=\begin{cases} \min \left\{s\ge 1 ~|~ s\binom{n+d-1}{n-1}-\binom {s}{2}\ge \binom{n+2d-1}{n-1}\right\}, & \text{ for } k=2;\\
\min \left\{s\ge 1 ~|~ s\binom{n+d-1}{n-1}\ge \binom{n+kd-1}{n-1}\right\}, & \text{ for } k\ge 3.
\end{cases}
\end{equation}
\end{conjecture}

Observe that, for $k\ge 3$, Conjecture~\ref{conj:main} claims that the na\"ive bound \eqref{ineq:lower bound} obtained by a parameter count is actually sharp, while, for $k = 2$, due to an additional group action there are many {\it defective} cases where the inequality is strict.

\begin{Rmk}
	Problems about additive decompositions including the above Waring problems can be usually rephrased geometrically in terms of {\em secant varieties}. In the case of $k$-Waring decompositions, we need to consider the {\em variety of powers} $V_{n,kd}^{(k)}$, i.e., the variety of $k$-th powers of forms of degree $d$ inside the (projective) space of forms of degree $kd$. The {\em $s$-th secant variety} $\sigma_s(V_{n,kd}^{(k)})$ is the Zariski  closure of the union of all linear spaces spanned by $s$-tuples of points lying on $V_{n,kd}^{(k)}$. In other words, it is the closure of the set of forms whose $k$-rank is at most $s$. Since the variety of powers is non-degenerate, i.e., it is not contained in any proper linear subspace of the space of forms of degree $kd$, the sequence of secant varieties  stratifies the latter space  and coincides with it for all sufficiently large $s$. Hence, the $k$-rank of a general form is the smallest value of $s$ for which the $s$-th secant variety of $V_{n,kd}^{(k)}$ coincides with the space of all forms of degree $kd$. Hence, Conjecture \ref{conj:main} can be rephrased as a conjecture about the dimensions of the secant varieties of $V_{n,kd}^{(k)}$. (We refer to \cite[Section 1.3.3]{On} for more details.) 
\end{Rmk}

\subsection{Maximal $k$-rank} 

A harder problem which is largely open even in the classical case of Waring decompositions, deals with the  computation of  the $k$-rank of an {\em arbitrary} complex-valued form of degree divisible by $k$. 

\begin{defi+}
Given a triple $(k,d,n)$,  denote by $\rk^{\max}_k(kd,n)$  the minimal number of terms such that {\em every}  form of degree $kd$  in $n+1$ variables can be represented as the sum of at most  $\rk^{\max}_k(kd,n)$ $k$-th powers of  forms of degree $d$. The number $\rk^{\max}_k(kd,n)$ is called the {\it maximal $k$-rank}.

(Similarly to the  above, we omit the subscript when considering the classical Waring rank, i.e., for $d = 1$.)
\end{defi+}
 
In \cite[Theorem 5.4]{Re1}, B. Reznick proved that the maximal Waring rank of binary forms of degree $k$ equals $k$. Moreover, the maximal value $k$ is attained exactly on the binary forms representable as $\ell_1 \ell_2^{k-1},$ where $\ell_1$ and $\ell_2$ are any two 
non-proportional  
linear binary forms. (Apparently these claims have  been known much earlier, but have never been carefully written down with a complete proof.)

\begin{problem}
Given a triple of positive integers $(k,d,n)$, calculate  $\rk_k^{\max}(kd,n).$
\end{problem}

At the moment, we have an explicit conjecture about the maximal $k$-rank only  in the case of binary forms. 

\begin{conjecture}\label{conj:maxrank}
For any positive integers $k,d$, the maximal $k$-rank $\rk^{\max}_k(kd,2)$ of binary forms equals $k$. Additionally, in the above notation, binary forms representable by $\ell_1 \ell_2^{kd-1}$, where $\ell_1$ and $\ell_2$ are non-proportional linear forms, have the latter maximal $k$-rank.
\end{conjecture}

Conjecture~\ref{conj:maxrank} is obvious for $k=2$ since, for any binary form $f$ of degree $2d$, we can write
\begin{equation} \label{eq:k2}
	f = g_1 g_2 = \left(\frac{1}{2} (g_1+g_2)\right)^2 + \left(\frac{i}{2} (g_1-g_2)\right)^2 \text{ with } g_1,g_2 \in S_d.
\end{equation}
The first non-trivial case is the one of binary sextics, i.e., $k=3, d=2$, which has been settled in \cite{LuOnReSh} where it has also  been  shown that the $4$-rank of $x_1x_2^7$ is equal to $4$.

\begin{Rmk} The best known general result about maximal ranks is due to G. Bleckherman and Z. Teitler, see \cite{BT} where they prove that the maximal rank is always at most twice as big as the generic rank. (This fact  is true  both for the classical ($d=1$) and for the higher ($d\geq 2$) Waring ranks.)

In the classical case of Waring ranks, this bound is (almost) sharp for binary forms, but in many other cases it is rather crude. At present, better bounds are known only in few special cases of low degrees \cite{BD13, Je14}. To the best of our knowledge, the exact values of the maximal Waring rank are only known for binary forms (classical, see \cite{Re1}), quadrics (classical), ternary cubics (see \cite{Seg42, LT10}), ternary quartics \cite{Kl99}, ternary quintics \cite{DeP15} and quaternary cubics \cite{Seg42}.
\end{Rmk}

%\begin{Rmk}{\rm The best known results about maximal rank in different situations are collected in \cite{BT}. In particular, they imply that $r_{max}(k,d,n)\le 2 r_{gen}(k,d,n)$.  In the case of binary forms, the latter inequality gives $$r_{max}(k,d,1)\le  2\left\lceil \frac{kd+1}{d+1}\right\rceil>k.$$
%}
%\end{Rmk}

\subsection{The $k$-rank of monomials}
Let $m = x_1^{a_1} \cdots  x_n^{a_n}$ be a monomial 
 with $0 < a_1 \leq a_2 \leq \cdots \leq a_n.$ It has been shown in \cite{ca-ch-ge} that the classical Waring rank of $m$ is equal to $\frac{1}{(a_1+1)}\prod_{i=1,\ldots,n} (a_i+1)$.

\medskip
Later E. Carlini and A.O. settled the  case of the $2$-rank, see \cite{ca-on}. Namely,  if $m$ is a monomial of degree $2d$, then we can write $m=m_1 m_2$, where $m_1$ and $m_2$ are monomials of degree $d$. From  identity (\ref{eq:k2}), it follows that the $2$-rank of $m$ is at most two. On the other hand, $m$ has rank  one exactly when we can choose $m_1 = m_2$, i.e., when  the power of each variable in  $m$ is even. 

While the cases $k=1$ and $k=2$ are solved, for $k\ge 3$, the question about the $k$-rank of monomials of degree $kd$,  is still open. At present, we are only aware of two general results in this direction. Namely, \cite{ca-on} contains the bound 
 $\rk_k(m) \leq 2^{k-1}$,  and recently, S.L., A.O., B.S., together with B. Reznick, have shown that $\rk_k(m) \leq k$ when $d \geq n(k-2)$, see \cite{LuOnReSh}. 
Thus, for fixed $k$ and $n$, all but a finite number of monomials of degree divisible by $k$ have $k$-rank less than $k$.
%At present we are unaware of any other general results on $k$-rank of monomials. % in this direction.

\begin{problem} \label{prob:monrk}
Given $k \geq 3$ and a monomial $m$ of degree $kd$, determine the monomial $k$-rank $\rk_k(m)$. 

 \end{problem}
%At present we are unaware of any other results in this direction. Even the case of binary forms is still open.) 
In the case of binary forms, a  bit more is currently known which motivates the following question.

\begin{problem}\label{prob:C}
Given $k \geq 3$ and a monomial $x^a y^b$ of degree $a+b=kd$, it is known that $\rk_k(x^{a}y^{b})  \leq \max(s,t)+1$, where $s$ and $t$ are the remainders of the division of  $a$ and $b$  by $k$, see  \cite{ca-on}. Is it true that the latter inequality is, in fact, an equality?  
\end{problem}
%Problem~\ref{prob:C} is a special case of the following more general problem about the $k$-ranks of monomials.

\subsection{Degree of the Waring map} 

Here again, we  concentrate on the case of binary forms (i.e., $n=2$). As we mentioned above, in this case, it is proven that 
$$\rk_k^\circ(kd,2)= \left \lceil \frac{\dim S_{kd}}{\dim S_{d}}\right \rceil= \left \lceil\frac{kd+1}{d+1}\right \rceil.$$

\begin{defi+}
We say that a pair $(k,d)$ is {\it perfect} if $\frac{kd+1}{d+1}$ is an integer. 
\end{defi+}

All perfect pairs are easy to describe. 

\begin{Lemma}\label{lm:perfect} The set of all pairs $(k,d)$ for which 
$\frac{kd+1}{d+1}\in\bN$ splits 
into the disjoint sequences 
$E_j := \{(jd+j+1,d) ~|~ d = 1,2,\ldots \}$. In each $E_j$, the corresponding quotient equals $jd+1$. 
\end{Lemma}

%For perfect pairs $(k,d),$ the following question is of substantial interest and difficulty. 

Given a perfect pair $(k,d)$, set $s:=\frac{kd+1}{d+1}$. Consider the map
$$W_{k,d}: S_d \times \ldots \times S_d \to S_{kd}, ~~(g_1,\ldots,g_s) \mapsto g_1^k + \ldots + g_s^k.$$
Let $\widetilde{W}_{k,d}$ be the same map, but defined up to a permutation of the $g_i$'s. We call it the {\it Waring map}. By \cite[Theorem 2.3]{LuOnReSh}, $\widetilde{W}_{k,d}$ is a generically finite map of complex linear spaces of the same dimension. By definition, its {\it degree} is the cardinality of the inverse image of a generic form in $S_{kd}$.

\begin{problem}\label{pr:deg}
Calculate the degree of $\widetilde{W}_{k,d}$ for  perfect pairs $(k,d)$.
\end{problem}

For the classical Waring decomposition ($d = 1$), we have a perfect pair if and only if $k$ is odd. From Sylvester's Theorem, we know that in this case the degree of the Waring map is $1$, i.e., the general binary form of odd degree has a {\it unique} Waring decomposition, up to a permutation of its summands.

\begin{Rmk}
	For the case of the classical Waring decomposition, the latter problem has  also been considered  in the case of more  variables. In modern terminology, the cases where the general form of a given degree has a unique decomposition up to a permutation of the summands are called {\it identifiable}. Besides the case of binary forms of odd degree, some other identifiable cases are classically known. These  are the quaternary cubics (Sylvester's Pentahedral Theorem \cite{Syl2}) and the ternary quintics \cite{Hilb, Pal03, Ri04, MM13}. Recently, F.~Galuppi and M.~Mella proved that these are the only possible identifiable cases, \cite{GM16}.
\end{Rmk}

\begin{Rmk}
	Problems dealing with additive decompositions of homogeneous polynomials similar to those we  consider in this section, have a very long story  going back to J.~J.~Sylvester and the Italian school of algebraic geometry of the late 19-th century. In the last decades, these problems received renewed attention due to their potential applications. Namely, homogeneous polynomials can be naturally identified with {\it symmetric tensors} and in several applied branches of science where such tensors are used, for example, to encode multidimensional data sets, {\it additive decompositions of tensors} play a crucial role as an efficient way to code those. 
		We refer to \cite{Lan} for an extensive exposition of these connections. 
\end{Rmk}

%Given a perfect pair $(k,d)$, set $s:=\frac{kd+1}{d+1}$. Consider the Waring map
%$$W_{k,d}: \oplus_{j=1}^s S^{d}_j \to S^{kd},$$
%sending an $s$-tuple of binary forms of degree $d$ to the sum of their $k$-th powers. Here $S^m$ stands for the linear space of binary forms of degree $m$. By the above result on the generic rank, $W_{k,d}$ is a generically finite map of complex linear spaces of the same dimension. Call be its {\it degree} the cardinality of the inverse image of a generic form in $S^{kd}$. 

%\begin{problem}\label{pr:deg}
%Calculate the degree of $W_{k,d}$ for  perfect pairs $(k,d)$.
%\end{problem}

%To the best of our knowledge, one of the very few situations when the answer to Problem~\ref{pr:deg} is known is Theorem~\ref{th:Sylv} (i) in which case such a representation is unique up to the order of summands implying that the degree of the corresponding Waring map equals $s!$. The major difficulty in Problem~\ref{pr:deg} is related to the fact that the map $W_{k,d}$ never projectivizes, i.e., one should deal with the case of excess intersection, see \cite {Fu}, Ch. 10.  

\section{Ideals of generic forms}
Let $I$ be a homogeneous ideal in $S$, i.e., an ideal generated by homogeneous polynomials. The ideal $I$ and the quotient algebra $R = S/I$ inherit the {\it grading} of the polynomial ring. 
%In particular, we have $I_i = I \cap S_i$ and $R_i = S_i/I_i$.%, for any $i \in \bN$. 
\begin{defi+}
	Given a homogeneous ideal $I \subset S$, we call the function 	$$
		%i \mapsto 
		\HF_R(i) := \dim_{\bC} R_i = \dim_{\bC} S_i - \dim_\bC I_i %,\quad \text{ for any } i \in \bN.
	$$
	the {\it Hilbert function} of $R.$
	The power series
	$$
		\HS_R(i) := \sum_{i\in\bN} \HF_R(i) t^i \in \bC[[t]]
	$$
	is called the {\em Hilbert series} of $R$.
\end{defi+}

Let $I$ be a homogeneous ideal generated by  forms $f_1, \dots, f_r$ of degrees $d_1, \dots, d_r$, respectively.
 It was shown in \cite{fr-lo} that, for fixed parameters $(n,d_1,\dots, d_r)$, there exists only a finite number of possible Hilbert series for $S/I$, and that
there is a Zariski open subset in the space of coefficients of the $f_i$'s
on which the Hilbert series of $S/I$ is one and the same and, in the appropriate sense, it is minimal among all possible Hilbert series, see below. We call algebras with this Hilbert series {\it generic}.
There is a longstanding conjecture about this minimal Hilbert series formulated by the first author, see \cite{fr}. 

\begin{conjecture}[Fr\"oberg's Conjecture, 1985]\label{conj:fr}
Let $f_1, \dots, f_r$ be generic forms of degrees  
$d_1, \dots, d_r$, respectively. Then the Hilbert series of the quotient algebra $R = S/(f_1,\ldots,f_r)$ is given by
\begin{equation}\label{eq:RALF}
\HF_R(t)=\left[\frac{\prod_{i=1}^r(1-t^{d_i})}{(1-t)^n}\right]_+.
\end{equation}
Here 
$[\sum_{i\ge0}a_iz^i]_+:=\sum_{i\ge0}b_iz^i$,
with $b_i=a_i$ if $a_j\ge0$ for all $j\le i$ and $b_i=0$. In other words,   $[\sum_{i\ge0}a_iz^i]_+$ is the truncation of a power series at its first non-positive coefficient.
\end{conjecture}

Conjecture~\ref{conj:fr} has been 
proven in the following cases: for $r\le n$ (easy exercise, since in this case $I$ is  a complete intersection); for $n\le2$, \cite{fr}; for $n=3$,  
\cite{an3}, for 
$r=n+1$,  which follows from \cite{st}. Additionally, in  \cite{ho-la} it has been proven that \eqref{eq:RALF} 
is correct in the first nontrivial degree $\min_{i=1}^r(d_i+1)$.  There are also other special results in the case  $d_1=\dots =d_r$, see  %\cite{ba-on},
 \cite{ni, au, fr-ho, mi-mi, ne}.
We should also mention that  \cite{fr-lu} contains a survey of the existing results on the generic series for various algebras and also it studies  the (opposite) problem of finding the maximal Hilbert series for fixed parameters $(n,d_1,\ldots,d_r)$.

\smallskip
%Although Conjecture~\ref{conj:fr} is not settled at present, there is a large number of related questions of interest. 

It is known that the actual Hilbert series of the quotient ring of any ideal with the same numerical parameters is lexicographically larger than or equal to the conjectured one. This fact implies that if for a  given discrete data $(n,  d_1, \ldots, d_r)$,  one finds just a single example of an algebra with the Hilbert series as in \eqref{eq:RALF}, then Conjecture~\ref{conj:fr} is settled in this case. 

\medskip
Although  algebras with the minimal Hilbert series constitute a Zariski open set, they are  hard to find constructively. We are only aware of two explicit constructions giving the minimal series  in the special case $r=n+1$, namely  R.~Stanley's choice $x_1^{d_1}, \ldots,x_n^{d_n}, (x_1+\cdots + x_n)^{d_{n+1}}$, and C. Gottlieb's choice  $x_1^{d_1}, \ldots,x_n^{d_n}, h_{d_{n+1}}$, where $h_{d}$ denotes the complete homogeneous symmetric polynomial of degree $d$, (private communication). To the best of our knowledge,  already in the next  case $r=n+2$ there is no  concrete guess about how to construct a similar example. There is however a substantial computer-based evidence pointing towards the possibility of replacing generic forms of degree $d$ by a product of generic forms of much smaller degrees. We present  some problems and conjectures related to such pseudo-concrete constructions below.

%\Samuel{Text above "final". Text below will be adjusted wrt Alessandros comments.}

\subsection{Hilbert series of generic power ideals.} Differently from the situation occurring in  R.~Stanley's result, if we consider ideals generated by more than $n+1$ powers of generic linear forms,  there are known examples of $(n,d_1,...,d_r)$ for which algebras generated by powers of generic linear forms  fail to have the Hilbert series as in \eqref{eq:RALF}.

\medskip
Recall that ideals generated by powers of linear forms are usually called {\it power ideals}. Due to their appearance in several areas of algebraic geometry, commutative algebra and combinatorics, they have been  studied more thoroughly. In the next section, we will discuss their relation with the so-called {\it fat points}. (For a more extensive survey of power ideals,  we refer to a nice paper by F. Ardila and A. Postnikov \cite{AP10}.) 

\medskip
Studying Hilbert functions of generic power ideals, A.~Iarrobino formulated the following conjecture, usually referred to as the Fr\"oberg-Iarrobino Conjecture, see \cite{ia, Ch}.

\begin{conjecture}[Fr\"oberg-Iarrobino Conjecture]\label{conj:fr-ia}
Given generic linear forms $\ell_1, \ldots, \ell_r$   and  a positive integer $d$, let $I$ be the power ideal generated by $\ell_1^d,\ldots,\ell_r^d$. Then the Hilbert function of $R = S/I$ is as in \eqref{eq:RALF}, except for the cases $(n,r) = (3,7), (3,8), (4,9), (5,14)$ and possibly for $r = n+2$ and $r = n+3$.
\end{conjecture}
This conjecture is still largely open. In \cite{fr-ho} R.F. and J. Hollman checked it for low degrees and low number of variables using the first version of the  software package {\it Macaulay2}. In the last decades, some progress has been made in reformulation of Conjecture~\ref{conj:fr-ia} in terms of the ideals of {\it fat points} and {\it linear systems}. We will return to this topic in the next section.

\subsection{Hilbert series of other classes of ideals}
Computer experiments suggest that in order to always generically get the Hilbert function as in \eqref{eq:RALF} we need to replace power ideals by slightly less special ideals.

\medskip
For example, given a partition $\mu = (\mu_1,\ldots,\mu_k) \vdash d$, we call  by a {\it $\mu$-power ideal}  an ideal generated by forms of the type $({\bf l}_1^\mu,\ldots,{\bf l}_r^\mu)$, where ${\bf l}_i^\mu=l_{i,1}^{\mu_1}\cdots l_{i,k}^{\mu_k}$ and $l_{i,j}$'s are distinct linear forms.
\begin{problem}\label{problem: F}
For $\mu \neq (d)$, does a generic $\mu$-power ideal have the same Hilbert function as in \eqref{eq:RALF}?
\end{problem}

Performed computer experiments  suggest a positive answer to the latter problem. L. Nicklasson has also conjectured that ideals generated  by powers of generic forms of degree $\geq 2$ have the Hilbert series as in \eqref{eq:RALF}.

%\begin{problem} \Samuel{Can be merged with problem 8}
%Let $f_1,\ldots,f_r$ be generic forms in $k[x_1,\ldots,x_n]$, $\deg(f_i)=d_i$. How many $f_i$ can be exchanged
%to powers of linear forms without changing the Hilbert series for the ideal they generate. \end{problem}

\begin{conjecture}[\cite{ni}] \label{conj:nic}
For generic forms $g_1,\ldots,g_r$   of degree $d>1$,  the ideal $(g_1^k,\ldots,g_r^k)$
has the same Hilbert series as the one generated by $r$ generic forms of degree $dk$. \end{conjecture}

 It was observed in \cite[Theorem A.3]{LuOnReSh} that Conjecture \ref{conj:nic} implies Conjecture \ref{conj:main}, connecting the two first sections of the present paper. It was also shown that Conjecture \ref{conj:nic} holds in the case of binary form by specializing the $g_i$'s to be $d$-th powers of linear forms and applying the fact that generic power ideals in two variables have the generic Hilbert series \cite{GeSh}. The same idea gives a positive answer to Problem \ref{problem: F}  in the case of binary forms, by specializing $l_{i,1} = \ldots = l_{i,k}$, for  $i = 1,\ldots,r$.

\subsection{Lefschetz properties of graded algebras}
We say that a graded algebra $A$ has the {\it weak Lefschetz property} (WLP) (respectively, the {\it strong Lefschetz property} (SLP)) if the multiplication map $\times l:A_i\rightarrow A_{i+1}$ (respectively,  $\times l^k:A_i\rightarrow A_{i+k}$) has the maximal rank, i.e., it is either injective or surjective, for a generic linear form $l$ and all $i$ (resp., for all $i$ and $k$). 
(For more references and open problems about the Lefschetz properties, see \cite{MiNa}.) % by J. Migliore and U. Nagel.

\begin{problem} It has been conjectured that each complete intersection $R=S/(f_1,\ldots,f_n)$ satisfies the WLP and also the SLP, see \cite{HMNW}. Does the same hold for $R=S/(f_1,\ldots,f_r)$, with $f_1,\ldots,f_r$ being generic forms, and $r > n$?
\end{problem}

It follows from \cite{St} that monomial complete intersections satisfy the SLP. In \cite{BFL}
the following situation has been studied. For  the ring $T_{n,d,k}=S/(x_1^d,\ldots,x_n^d)^k$, it is shown that for $k\ge d^{n-2}$, $n\ge3$, $(n,d)\ne(3,2)$,  $T_{n,d,k}$
fails the  WLP. For $n=3$, there is an explicit  conjecture when the WLP holds.  Additionally,  there is some information about  $n>3$.

\begin{problem} When are the WLP and the SLP true for $T_{n,d,k}$?
\end{problem}

We now introduce the concept of the \emph{$\mu$-Lefschetz properties}. 
%\begin{defi+}
Let $\mu=(\mu_1,\ldots,\mu_k)$ be a partition of $d$, i.e.,  $\sum_{i=1}^k\mu_i=d$. We say that an algebra has the {\it $\mu$-Lefschetz property} if $\times{\bf l}^\mu:A_i\rightarrow A_{i+d}$ has maximal rank for all $i$, where ${\bf l}^\mu=l_1^{\mu_1}\cdots l_k^{\mu_k}$, and $l_i$'s are 
generic linear forms.
%\end{defi+}
 
\begin{problem}
For $R= S/(f_1,\ldots,f_r)$, where $f_1,\ldots,f_r$ are generic forms, does $R$ satisfy the $\mu$-Lefschetz property for all partitions $\mu$?
 \end{problem}

\section{Symbolic powers}

%\begin{defi+}
For a prime ideal $\wp$ in a Noetherian ring $R$, define  its {\it $m$-th symbolic power} $\wp^{(m)}$ as  $$\wp^{(m)}=\wp^mR_\wp\cap R.$$ 
It is the $\wp$-primary component of $\wp^m$. For a general ideal $I$ in $R$, its {\it $m$-th symbolic power} is defined as
$I^{(m)}=\cap_{\wp\in{\rm Ass}(I)}(I^mR_\wp\cap R)$. 

\subsection{Hilbert functions of fat points.}
Let $I_X$ be the ideal in $\bC[x_1,\ldots,x_n]$ defining a scheme of reduced points $X = P_1 + \ldots + P_s$ in $\bP^{n-1}$, say $I_X = \wp_1 \cap \ldots \cap \wp_s$ where $\wp_i$ is the prime ideal defining the point $P_i$. Then, the $m$-th symbolic power $I^{(m)}$ is the ideal $I_X^{(m)} = \wp_1^m \cap \ldots \cap \wp_s^m$ which defines the scheme of {\it fat points} $X = mP_1 + \ldots + mP_s$. 

Ideals of $0$-dimensional schemes are  classical objects of study since the beginning of the last century.  Their Hilbert functions are of particular interest. Study of these ideals and calculation of their Hilbert functions can be often  related to the so-called {\it polynomial interpolation problem}. Indeed, the homogeneous part of degree $d$ of the ideal $I_X$ is the space of hypersurfaces of degree $d$ in $\bP^{n-1}$ passing through the $P_i$'s up to order $m-1$, i.e., the space of polynomials of degree $d$ whose partial differentials up to order $m-1$ vanish at every $P_i$.

It is well-known that the Hilbert function of $0$-dimensional schemes is strictly increasing until it reaches the {\it multiplicity} of the scheme, see \cite[Theorem 1.69]{IK06}. Hence, since the degree of a $m$-fat point in $\bP^{n-1}$ is ${n-1+m-1 \choose n-1}$, the expected Hilbert function is $$\HF_{S/I_X}(d) = \min\left\{{n-1+d \choose n-1}, s{n-1+m-1 \choose n-1}\right\}.$$ In the case of simple generic points, i.e., for $m=1$, it is known that  the actual Hilbert function is as expected.

In the case of double points ($m=2$), counterexamples were known since the end of the 19-th century. In 1995, after a series of important papers, J. Alexander and A. Hirschowitz proved that the classically known examples were the only counterexamples. For higher multiplicity, very little is known at present. In the case of projective plane, a series of equivalent conjectures have been given by  B. Segre \cite{Seg61}, B. Harbourne \cite{Har86}, A. Gimigliano \cite{Gim87} and A. Hirschowitz \cite{Hir89}. These are known as the {\it SHGH-Conjecture}, see \cite{Har00} for a survey of this topic.

\smallskip
{\it Apolarity Theory} is a very useful tool in studying of the ideals of fat points and it is  connecting  all the algebraic stories we have told above. In particular, the following lemma is crucial. 
(We refer to \cite{IK06} and \cite{Ger96} for an  extensive description of this issue.)

\begin{Lemma}[Apolarity Lemma]
	Let $X = P_1+\ldots+P_s$ be a scheme of reduced points in $\bP^{n-1}$ and let $L_{1},\ldots, L_{s}$ be linear forms in $\bC[x_0,\ldots,x_n]$ such that, for any $i$, the coordinates of $P_i$ are the coefficients of $L_i$. Then, for every $m\geq d$,
	$$
		\HF_{S/I^{(m)}_X}(d) = \dim_{\bC}[(L_1^{m-d+1},\ldots,L_s^{m-d+1})]_d.
	$$
\end{Lemma}

Using this statement we obtain that calculation of  the Hilbert function of a scheme of fat points is equivalent to the calculation of  the Hilbert function of the corresponding power ideal. In particular, Fr\"oberg-Iarrobino conjecture (Conjecture \ref{conj:fr-ia}) can be rephrased as a conjecture about the Hilbert function of ideals of generic fat points.

\medskip
Recently R.F. raised the question about what happens in case of the ideals of generic fat points in a multi-graded space. A point in multi-projective space $P \in \bP^{n_1-1}\times\ldots\times \bP^{n_t-1}$ is defined by a prime ideal $\wp$ in the multi-graded polynomial ring $S = \bC[x_{1,1},\ldots,x_{1,n_1};\ldots;x_{t,1},\ldots,x_{t,n_t}] = \bigoplus_{I \subset \bN^t} S_I$, where $S_I$ is the vector space of multi-graded polynomials of multi-degree $I = (i_1,\ldots,i_t) \in \bN^t$. A scheme of fat points $X = mP_1+\ldots+mP_s$ is the scheme associated with the multi-graded ideal $\wp_1^m\cap\ldots\cap\wp_s^m$.

\begin{problem}
	Given a scheme of generic fat points $X \subset \bP^{n_1-1}\times\ldots\times \bP^{n_t-1}$, what is the multi-graded Hilbert function $\HF_{S/I_X}(I)$, for $I \in \bN^t$?
\end{problem}

This question was first considered by M. V. Catalisano, A. V. Geramita and A. Gimigliano who solved it in the case of double points, i.e., for $m=2$ in $\bP^1 \times \ldots \times \bP^1$. Recently, A.O. jointly with E. Carlini and M. V. Catalisano    resolved the case of triple points ($m=3$) in $\bP^1 \times \bP^1$ and computed the Hilbert function for an arbitrary multiplicity except for a finite region in the space of multi-indices, see \cite{CCO17}. %\Alessandro{As far as we know, these are the only known cases.}

\subsection{Symbolic powers vs. ordinary powers.}
As we mentioned above, if $I$ is the ideal defining  a set $X$ of points, the $m$-th symbolic power of $I$ is the ideal of polynomials vanishing up to order $m-1$ at all points in $X$ or, in other words, the space of hypersurfaces which are singular at all points in $X$ up to order $m-1$. For this reason,  symbolic powers are interesting from a geometrical point of view, but they are more difficult to study  compared to the usual powers which carry  less geometrical information.  Hence, it is important to find relations between them.  Observe that the inclusion $I^m \subset I^{(m)}$ is trivial.

\smallskip
{\it Containment problems} between the ordinary and the symbolic powers of ideals of points have been  studied in substantial details. One particularly interesting question is  to understand for which pairs of positive integers $(m,r)$,  $I^{(m)} \subset I^r$. A very important result in this direction  is the fact that, for any ideal $I$ of reduced points in $\bP^n$ and any $r > 1$, we have $I^{(nr)}\subset I^r$. This statement was proven in \cite{ELS01} by L. Ein, R. Lazersfeld and K. Smith for characteristic $0$  and by M. Hochster and C. Huneke in positive characteristic, see \cite{HH02}. At present, the important question is whether the bound in the latter statement is sharp. In \cite{DSTG13}, M. Dumicki, T. Szemberg and H. Tutaj-Gasi\'nska provided the first example of a configuration of points such that $I^{(3)} \not\subset I^2$. (We refer to \cite{SS17} for a complete account on this topic.) 

\medskip
In the recent paper \cite{GGSVT16}, F. Galetto, A. V. Geramita, Y. S. Shin and A. Van Tuyl defined the {\it $m$-th symbolic defect of an ideal} as the number of minimal generators of the quotient ideal $I^{(m)} / I^m$. If $I$ defines a set of general points in projective space, it was already known that $I^{(m)} = I^m$ if and only if  $I$ is a complete intersection. Additionally, in \cite[Theorem 6.3]{GGSVT16} the authors characterize all cases of $s$ points in $\bP^2$ having the $2$-nd symbolic defect equal to $1$. These cases are exactly $s = 3,5,7,8$.

\begin{problem} For the ideal $I$ of $s$ general points in $\mathbb P^{n-1}$, what is the difference between the Hilbert series of the $m$-th symbolic power and the $m$-th ordinary power? 
\end{problem}

%
%the difference between the Hilbert 
%series of $I^{(2)}$ and $I^2$ for $s=1,2,\ldots,10$ is $0,0,t^3,0,t^5,3t^5,t^6,t^6,3t^7,6t^7$.
%Could this be explained? (For $s=1,2,4$ we have complete intersections. For e.g. 18
%points the difference is $3t^{11}+6t^{10}+t^9$.)\end{problem}
%%\begin{problem} Can one say something about the primary decomposition for ideals of generic points in $\mathbb P^2$?
%\end{problem}

\section{Miscellanea}  
\subsection{Hilbert series of numerical semigroup rings}
%\subsection{}
Let $\mathcal{S}=\langle s_1,\ldots,s_k\rangle$ be a numerical semigroup, i.e. $\mathcal{S}$ consists of all linear combinations with non-negative integer coefficients of the positive integers
$s_i$, and let $k[x^{s_1},\ldots,x^{s_k}]=k[\mathcal{S}]$ be the semigroup ring. The Hilbert series of $k[\mathcal{S}]$ is of the form $p(t)/q(t)$, where $p,q$ are polynomials with integer coefficients.
A semigroup is called \emph {cyclotomic} if the polynomial $p(t)$ has all its roots in the unit circle (which in fact, implies that they lie on the unit circle). (Detailed information about numerical semigroups can be found in  e.g., \cite{Ci}.) 

\begin{conjecture} $\mathcal{S}$ is cyclotomic if and only if $k[\mathcal{S}]$ is a complete intersection.
\end{conjecture}

\subsection{Non-negative forms}  The next circle of problems is related to the celebrated article  \cite{Hi} of D.~Hilbert and to a number of results formulated in  \cite{CLR}.   

Denote by  $P_{n,m}$ the set of all non-negative real forms, i.e., real homogeneous polynomials  of (an even) degree $m$ in $n$ variables which never attain negative values; denote by $\Sigma_{n.m}\subseteq P_{n,m}$ the subset of non-negative forms which can be represented as sums of squares of real forms of degree $\frac{n}{2}$.  
(In \cite{Hi} D.~Hilbert proved that $\Delta_{n,m}=P_{n,m}\setminus\Sigma_{n,m}$ is non-empty unless the pair $(n,m)$ is of the form $(n,2)$, $(2,m)$  or $(4,3)$.)  Finally, if $\Z(p)$ stands for the real zero locus  of a real form $p$,  denote by $B_{n,m}$ (resp. $B'_{n,m}$) the supremum of $|\Z(p)|$ over $p\in P_{n,m}$ such that $|\Z(p)|<\infty$  
(resp. over $p\in\Sigma_{n,m}$ such that $|\Z(p)|<\infty$). In other words, $B(n,m)$ is the supremum of the number of zeros of non-degenerate forms under the assumption that all these roots are isolated (and similarly for 
$B'_{n,m}$).  Obviously, $B'_{n,m}<B_{n,m}$.

\medskip 
The following basic question  was posed in \cite{CLR}. 

\begin{problem} Are  $B_{n,m}$ and $B'_{n,m}$ finite for any pair $(n,m)$ with even $n$? \end{problem}

In \cite{CLR} it was shown that the answer to this problem is  positive  for $m=2,3$ and for the pair $(4,4)$. Relatively recently,  in \cite{Stu} the following upper bound for $B_{n,m}$ was established
$$B_{n,m}\le 2\frac{(m-1)^{n+1}-1}{m-2}.$$
However this bound can not be sharp, as shown in \cite{Ko}.  In case of $B'_{n,m}$, the following guess seems  quite plausible and is proven for $m=3$.

\begin{conjecture}  For any given pair $(n,m)$ with even $n$, $B'_{n,m}=\left(\frac{n}{2}\right)^{m-1}$.

\end{conjecture}

For $B_{n,m}$, no similar guess is known, but some intriguing information  is available  in the case $m=3$, see \cite{CLR}. The following problem is related to the classical Petrovski-Oleinik upper bound on the number of real ovals of real plane algebraic curves.

\begin{problem} Determine $\lim_{n\to\infty}\frac{B_{n,3}}{n^2}$. \end{problem}

The latter limit exists and lies in the interval  $\left[\frac{5}{18},\frac{1}{2}\right]$, see  \cite{CLR}. 

%\begin{problem} Is $B_{4,4}=10$ or 11?\end{problem}

%\begin{problem} Is $B'_{4,4}=8$ or 9 or 10?\end{problem}

\subsection{Polynomial generation}
Let $p$ be a prime number and let $\mathbb{F}_p$ denote the field with $p$ elements. Consider the two maps 
$$\phi: \mathbb{F}_p[x_1,\ldots,x_n] \to \mathbb{F}_p[x_1,\ldots,x_n], f \mapsto \sum_{a \in Z(f)} x^a,$$ 
$$\psi: \mathbb{F}_p[x_1,\ldots,x_n] \to \mathbb{F}_p[x_1,\ldots,x_n], f \mapsto \sum_{a \in \mathbb{F}_p^n} f(a) x^a.$$
Here  $x^a := x_1^{a_1} \cdots x_n^{a_n}$, where each $a_i$ is regarded as an integer, and $Z(f)$ is the zero locus of $f$ in $F_p^n$, i.e.,  $Z(f) := \{a \in \mathbb{F}_p^n \,|\, f(a) = 0\}$.
When $p = 2$, then $\phi$ is a bijective map on the vector space of polynomials of degree at most one in each variable, and $\phi^4(f) = f$, see \cite{lu}.
The map $\psi$, suggested by M. Boij, is a linear bijective map on the vector space of polynomials with degree at most $p-1$ in each variable, and when $p = 2$, these two maps are closely related in the sense that
$\phi(f) = \psi(f) + \sum_{a \in \mathbb{F}_2^n} x^a.$

Consider now the case $n = 1$ and $p > 2$. The map $\phi$ is no longer a bijection, but the sequence $\phi(f),\phi^2(f),\ldots$ will eventually become periodic. It is an easy exercise to show that $0 \mapsto 1 + x + \cdots + x^{p-1} \mapsto x \mapsto 1 \mapsto 0$. When $p \leq 17$, this is the only period, i.e., $\phi(f)^{d(f)} = 0$ for some $d(f)$. For $p=71$,  we have found a period of length two; 
$1 +x^{63}
\mapsto x^{23} + x^{26} + x^{34} + x^{39} + x^{41} + x^{51} + x^{70} \mapsto 1 + x^{63}.$ One can show that the length of the period is always an even number, but it is not clear which even numbers that can occur as lengths of periods.

\begin{problem}
For $n=1$ and given $p$, what are the (lengths of the) possible periods of $\phi$? 
\end{problem}

Let us now turn to the map $\psi$ and the case $n=1$. For $p =3$,  $\psi^8(f) = f$ for all polynomials $f$ in $\mathbb{F}_3[x]$ of degree at most two.  For $p = 5$, the least $i$ such that $\psi^{i} = {\rm Id}$ on the space of polynomials of degree at most four, is equal to $124$. For $p = 7$, the corresponding number is $1368$.

\begin{problem}
For $n = 1$ and given $p$, find the minimal positive integer  $i$ such that $\psi^i$ is the identity map on the space of polynomials of degree at most $p-1$.
\end{problem}

\subsection{Exterior algebras}
Let $f$ be a generic form of even degree in the exterior algebra $E$ over $\mathbb{C}$ with $n$ generators. Moreno and Snellman showed that the Hilbert series of $E/(f)$ is equal to the expected series $[(1+t)^d(1-t^d)]_+$, see \cite{M-S}. When the degree of $f$ is odd, we have $(f)\subseteq{\rm Ann}{f}$. This annihilator ideal shows an unexpected behaviour. The most striking case is when 
$(n,d) = (9,3)$. It turns out that $\dim_{\mathbb{C}}({\rm Ann}{f})_3 = 4$, see \cite{lu-ni}. At the same time a naive guess is that $({\rm Ann}{f})_3$ is spanned by $f$ only.
Additionally, computer experiments suggest that $(f)$ and ${\rm Ann}{f}$ agree in low degrees.

%The problem suggested below comes from \cite{lu-ni}.
\begin{problem} Let $f$ be a form of odd degree $d$ in $E$.  Is it true that $({\rm Ann}(f))_i = (f)_i$, for $i < (n-d)/2$? \end{problem}

\medskip
We finish our list of problems  with the following conjecture stated in  \cite{cr-lu-ne}, which connects the question about the Hilbert series of generic forms in the exterior algebra with the Hilbert series of power ideals in the commutative setting.
\begin{conjecture}
Let $f$ and $g$ be generic quadratic forms in $E$ and let $\ell_1$ and $\ell_2$ be two generic linear forms in $S$. Then the Hilbert series of $E/(f,g)$ is equal to the Hilbert series of $S/(x_1^2,\ldots,x_n^2,\ell_1^2,\ell_2^2)$ and is given by 
$1 + a(n,1) t + a(n,2) t^2 + \cdots + a(n,s)t^s + \cdots$, where $a(n,s)$ is the number of lattice paths inside the rectangle $(n+2-2s)\times (n+2)$ starting from the bottom-left corner and ending at the top-right corner by using only moves of two types: either $(x,y)\rightarrow (x+1,y+1) \textrm{\ or\ } (x-1,y+1)$.

\end{conjecture}

%so the situation is more complicated. 

\medskip\noindent
{\bf Acknowledgements.}   The authors want to  thank all the participants of the problem-solving seminar at Stockholm University for their contributions and patience. The 
fourth author is sincerely grateful to Dr.~Kh.~Khozhasov for pointing out references \cite{Stu, Ko}.


\begin{thebibliography}{ Dillo 83}

\bibitem[AH95]{al-hi} Alexander, J., and Hirschowitz, A. {\it Polynomial interpolation in several variables}, J. Alg. Geom. {\bf 4}, 201--222 (1995).

\bibitem[An86]{an3} Anick, D. {\it Thin algebras of embedding dimension three}, J. Algebra {\bf 100},
235--259 (1986).

\bibitem[AP10]{AP10} Ardila, F., and Postnikov, A. {\it Combinatorics and geometry of power ideals}, Transactions of the American Mathematical Society {\bf 362}(8), 4357--4384, (2010).

\bibitem[Au95]{au} Aubry, M. {\it S\'erie de Hilbert d'une alg\`ebre de polyn\^omes quotient}, J. Algebra {\bf 176}, 392--416 (1995).

\bibitem[BD13]{BD13}
Ballico, E., and De Paris, A. {\it Generic power sum decompositions and bounds for the Waring rank}, Discrete \& Computational Geometry {\bf 57}(4), 896-914. (2017).

\bibitem[BT15]{BT} Blekherman, G., and Teitler, Z. {\it On maximum, typical, and generic ranks}, Mathematische Annalen, {\bf 362}(3--4),  1021--1031 (2015).

\bibitem[BFL18]{BFL} Boij, M., Fr\"oberg, R., and Lundqvist, S. {\it Powers of generic ideals
and the weak Lefschetz property for powers of monomial complete intersections}, J. Algebra
{\bf 495} 1--14 (2018).
%\bibitem[Br-Ot]{BO} Brambilla, M.~C. Ottaviani, G. {\it On the Alexander-Hirschowitz theorem}, J.Pure Appl. Algebra, {\bf 212}, 1229--1251  (2008).

\bibitem[CCG12]{ca-ch-ge}
 Carlini, E.,  Catalisano, M. V.,  and  Geramita, A. V. {\it The solution to the Waring problem
for monomials and the sum of coprime monomials}, J. Algebra {\bf 370}, 5--14, (2012).

\bibitem[CCO17]{CCO17}
Carlini, E., Catalisano M. V., and Oneto A. {\it On the Hilbert function of general fat points in $\mathbb {P}^ 1\times\mathbb {P}^ 1$}, arXiv preprint arXiv:1711.06193 (2017).

\bibitem[CO15]{ca-on}
Carlini, E., and  Oneto, A. {\it Monomials as sum of k-th powers of forms}, Comm. Algebra {\bf 43}, 650--658 (2015).

\bibitem[CS13]{Stu}  Cartwright, D., and  Sturmfels, B. {\it The number of eigenvalues of a tensor}, Linear Algebra and its Applications {\bf 438}, 942--952 (2013).

%\bibitem[CGGH+15]{many} Catalisano, M.~V., Geramita, A.~V., Gimigliano, A., Harbourne, B.,  Migliore, J., Nagel, U., Shin, Y.~S. {\it Secant varieties of reducible varieties in ${\mathbb P}^n$}, arXiv: 1502.00167 (2015).

\bibitem[Ch05]{Ch} Chandler, K.~{\it The geometric interpretation of Fr{\"o}berg--Iarrobino conjectures on infinitesimal neighbourhoods of points in projective space}, J. of Algebra {\bf 286}(2),  421--455, (2005).

\bibitem[CLR80]{CLR} Choi, M.-D.,  Lam, T.-Y., and Reznick, B.~{\it Real Zeros of Positive Semidefinite Forms. I}, Math. Z., {\bf 171}, 1--26 (1980).

\bibitem[CGM16] {Ci} Ciolan, A.,  Garcia-Sanchez, P., and Moree, P. {\it Cyclotomic numerical semigroups},  SIAM J. Discrete Math., {\bf 30}(2), 650--668 (2016). 

\bibitem[CLN]{cr-lu-ne} Crispin Qui{\~n}onez, V.,   Lundqvist, S., and  Nenashev, G. {\it On ideals generated by two generic quadratic forms in the Exterior algebra}, in preparation.

\bibitem[DeP15]{DeP15}
De Paris, A. {\it Every ternary quintic is a sum of ten fifth powers}, International Journal of Algebra and Computation (2015)

\bibitem[DSTG13]{DSTG13}
Dumnicki, M.,  Szemberg, T., and  Tutaj-Gasi\'nska H. {\it Counterexamples to the $^{(3)}\subset I^2$ containment}. J. Algebra {\bf 393}, 24--29 (2013).

\bibitem[ELS01]{ELS01}
 Ein, L.,  Lazarsfeld, R., and  Smith, K.  {\it Uniform bounds and symbolic powers on smooth varieties}. Invent. Math. {\bf 144} 241-- 252 (2001).

\bibitem[Fr85]{fr} Fr\"oberg, R. {\it An inequality for Hilbert series.} Math. Scand. {\bf 56}, 117--144 (1985).

\bibitem[FH94]{fr-ho} Hollman, J., and Fr\"oberg, R. {\it Hilbert series for ideals generated by generic forms}, J. Symb. Comp. {\bf 17}, 149--157 (1994).

\bibitem[FL90]{fr-lo} Fr\"oberg, R., and L\"ofwall, C. {\it On Hilbert series for commutative and noncommutative graded algebras}, J. Pure Appl. Algebra {\bf 76}, 33--38 (1990).

\bibitem[FL17]{fr-lu} Fr\"oberg, R., and Lundqvist, S. {\it Extremal Hilbert series.} arXiv: 1711.01232, (2017).

\bibitem[FOS12]{FOS} Fr\"oberg, R., Ottaviani, G., and Shapiro, B.~ {\it On the Waring problem for polynomial rings},  PNAS,  {\bf 109}(15), 5600--5602 (2012).

%\bibitem[Fu84]{Fu} Fulton, W. Intersection theory,  Ergebnisse der Mathematik und ihrer Grenzgebiete 3.  Band 2, Springer-Verlag, 1984, xi + 470 pp.

\bibitem[GGSVT16]{GGSVT16}
Galetto, F., Geramita, A.~V., Shin, Y.~S. and Van Tuyl, A. {\it The symbolic defect of an ideal.}  arXiv:1610.00176, (2016).

\bibitem[Ger96]{Ger96}
Geramita, A.~V. {\it Waring's Problem for Forms: inverse systems of fat points, secant varieties and Gorenstein algebras.} Queen's Papers in Pure and Applied Math. {\bf 105}(2),  1--129,  (1996).

\bibitem[Gim87]{Gim87}
 Gimigliano, A. {\it On Linear Systems of Plane Curves}, Ph. D. thesis, Queen's University, Kingston, Ontario (1987).

\bibitem[GM16]{GM16}
Galuppi, F., and  Mella, M. {\it Identifiability of homogeneous polynomials and Cremona Transformations},  arXiv:1606.06895 (2016).

\bibitem[GS98]{GeSh}
Geramita, A.~V., and Schenck, H.~ {\it Fat points, inverse systems, and piecewise polynomial functions}, Journal of Algebra {\bf 204}(1), 116--128 (1998).

\bibitem[Har86]{Har86}
 Harbourne, B. {\it The geometry of rational surfaces and Hilbert functions of points in the plane}, Can. Math. Soc. Conf. Proc. {\bf 6}, 95--111 (1986).

\bibitem[Har00]{Har00}
Harbourne, B. {\it Problems and progress: a survey on fat points in P2}, in Zero Dimensional Schemes and Applications, Naples (2000), 85--132; Queenâs Papers in Pure and Appl. Mat.,

\bibitem[HMNW03]{HMNW} Harima, T., Migliore J., Nagel, U., and Watanabe J. {\it The weak
and strong Lefschetz properties for artinian $K$-algebras}, J. Algebra {\bf 262}, 99--126 (2003).

\bibitem[Hil88]{Hi} D.~Hilbert,  {\it \"Uber die Darstellung definiter Formen als Summe von Formenquadraten,} Math. Ann. {\bf 32}(3), 342--350 (1888).

\bibitem[Hilb]{Hilb}
 Hilbert, D. {\it Letter adress\'ee \`a M. Hermite}, Gesam. Abh. vol II, 148--153.

\bibitem[Hir89]{Hir89}
 Hirschowitz, A. {\it Une conjecture pour la cohomologie des diviseurs sur les surfaces rationelles ge\'en\'eriques}, Journ. Reine Angew. Math. {\bf 397} 208--213 (1989).

\bibitem[HL87]{ho-la} Hochster, M. and Laksov, D. {\it The linear syzygies of generic forms}, Comm. 
Algebra {\bf 15}, 227--234 (1987).

\bibitem[HH02]{HH02}
 Hochster, M.,  and  Huneke, C. {\it Comparison of symbolic and ordinary powers of ideals.} Invent. Math. {\bf 147} 349 -- 369 (2002).
 
\bibitem[Ia97]{ia} Iarrobino, A. {\it Inverse system of a symbolic power III. Thin algebras and fat points}, Compositio Math. {\bf 108}, 319--356 (1997).

\bibitem[IK06]{IK06}
Iarrobino, A., and Kanev, V. {\it Power sums, Gorenstein algebras, and determinantal loci.} Springer (2006).

\bibitem[Kl99]{Kl99}
 Kleppe, J. {\it Representing a homogeneous polynomial as a sum of powers of linear forms,} Master Thesis, University of Oslo (1999).

\bibitem[Ko17]{Ko} Kozhasov, Kh., {\it On fully real eigenconfigurations of tensors}, https://arxiv.org/pdf/1707.04005.pdf. 

\bibitem[Je14]{Je14} 
Jelisiejew, J. {\it An upper bound for the Waring rank of a form,} Archiv der Mathematik {\bf 102}(4), 329--336 (2014).

%\bibitem[Ho-Ki-Va] {HKV} A.~Hosry, Y.~ Kim and J.~Validashti, {\it On the equality of ordinary and symbolic powers of ideals},   J. Commut. Algebra {\bf 4}(2), 281--292   (2012).
\bibitem[La12]{Lan}
Landsberg, J. M. {\it Tensors: geometry and applications}, Vol. 128. Providence, RI: American Mathematical Society, 2012.

\bibitem[Lu15]{lu} Lundqvist, S. {\it Boolean Ideals and their Varieties}, J. Pure Appl. Alg. {\bf 219}(5), 4521--4540 (2015).

\bibitem[LN]{lu-ni} Lundqvist, S., Nicklasson, L., {\it On generic principal ideals in the Exterior algebra},  in preparation.

\bibitem[LORS17]{LuOnReSh} Lundqvist, S., Oneto, A.,  Reznick, B., and Shapiro, B.  {\it On generic and maximal k-ranks of binary forms},  arXiv:1711.05014, submitted.

\bibitem[LT10]{LT10}
Landsberg, J. M, and Teitler, Z. {\it On the ranks and border ranks of symmetric tensors}, Foundations of Computational Mathematics {\bf 10}(3), 339--366 (2010).

\bibitem[MM13]{MM13}
 Massarenti, A.,  Mella, M. {\it Birational aspects of the geometry of varieties of sums of
powers}, Adv. Math. {\bf 243}, 187--202 (2013).

\bibitem[MM03]{mi-mi} Migliore, J., and Miro-Roig, R. M. {\it Ideals of general forms and the ubiquity of the weak Lefschetz property}, J.  Pure Appl. Algebra {\bf 102}, 79--107 (2003).

%\bibitem[MM02]{Mi-Mi2} Migliore, J. and Miro-Roig, R. M. {\it On the minimal free resolution of $n + 1$ general forms}, Trans. Amer. Math. Soc. {\bf 355},1--66 (2002).


\bibitem[MN03]{MiNa} Migliore, J., and Nagel, U. {\it The Weak and Strong Lefschetz properties for Artinian K-algebras}, J. Algebra, {\bf 262}, 99--126  (2003). 

%\bibitem{Mo}  J.~Molluzzo,  Monotonicity of quadrature formulae and polynominal representation [Doctoral thesis]. Yeshiva University; 1972.

\bibitem[MS02]{M-S} Moreno-Socias, G., and Snellman, J.~ {\it Some conjectures about the Hilbert series of generic ideals in the exterior algebra},
Homology Homotopy Appl. {\bf 4}, 409--426 (2002).

\bibitem[Ne17]{ne} Nenashev, G. {\it A note on Fr\"obergs's conjecture for forms of equal degree},
Comptes Rendus Mathematique, {\bf 355}(3), 272--276  (2017).

%\bibitem[NSS]{NeShSh} Nenashev, G., Shapiro, B., and Shapiro, M. {\it Secant degeneracy index of the standard strata in the space of binary forms}, Arnold Mathematical Journal, to appear.

\bibitem[Ni17]{ni} Nicklasson, L. {\it On the Hilbert series of ideals generated by generic forms},
Comm. Algebra, {\bf 45}(8), (2017).

\bibitem[On16]{On} Oneto, A.~ {\it Waring-type problems for polynomials}, Doctoral Thesis in Mathematics, Stockholm University, Stockholm, Sweden (2016).

\bibitem[Pal03]{Pal03}
 Palatini, F. {\it Sulla rappresentazione delle forme ternarie mediante la somma di potenze
di forme lineari}, Rom. Acc. L. Rend. {\bf 12}, 378--384 (1903).

\bibitem[Re13]{Re1} Reznick, B.~ {\it Some new canonical forms for polynomials}, Pacific J.Math. {\bf 266}(1)  185--220 (2013).

\bibitem[Ri04]{Ri04}
Richmond, H. W.  {\it On canonical forms}, Quart. J. Pure Appl. Math. {\bf 33}, 967--984 (1904).

\bibitem[Seg42]{Seg42}  Segre, B. {\it The Non-singular Cubic Surfaces}, Oxford University Press, Oxford, 1942.

\bibitem[Seg61]{Seg61}
 Segre, B. {\it Alcune questioni su insiemi finiti di punti in Geometria Algebrica}, Atti del Convegno Internaz. di Geom. Alg., Torino (1961).

\bibitem[St78]{st} Stanley, R. {\it Hilbert functions of graded algebras}, Adv. Math. {\bf 28} 57--83 (1978).

\bibitem[St80]{St} Stanley, R. {\it Weyl groups, the hard Lefschetz theorem, and the Sperner
property}, SIAM J. Algebraic Discrete methods {\bf 1}(2), 168--184 (1980).

\bibitem[Sy51]{Syl2}   Sylvester, J. J. {\it On a remarkable discovery in the theory of canonical forms and of hyperdeterminants,}
originally in Philosophical Magazine, vol. I, 1851; pp. 265--283 in Paper 41 in
Mathematical Papers, Vol. 1, Chelsea, New York, 1973. Originally published by Cambridge
University Press in 1904.

\bibitem[SS17]{SS17}
Szemberg, T., and Szpond J. {\it On the containment problem.} Rendiconti del Circolo Matematico di Palermo Series 2, {\bf 66}(2), 233--245  (2017).

\end{thebibliography}
\end{document}